        \documentstyle[amsmath,amsfonts,euscript]{article}
        \newcommand{\al}{\alpha}
        \newcommand{\gam}{\gamma}
        \newcommand{\del}{\delta}
        \newcommand{\eps}{\epsilon}
        \newcommand{\veps}{\varepsilon}
        \newcommand{\lam}{\lambda}
        
        \newcommand{\sig}{\sigma}
        
        \newcommand{\om}{\omega}
        
        \newcommand{\Gam}{{\it \Gamma}}
        \newcommand{\Lam}{{\it \Lambda}}
	\newcommand{\XI}{{\it \Xi}}
	\newcommand{\Om}{{\it \Omega}}

        \newcommand{\NN}{{\mathbb{N}}}
        \newcommand{\ZZ}{{\mathbb{Z}}}
        \newcommand{\RE}{{\mathbb{R}}}
        \newcommand{\CO}{{\mathbb{C}}}
	\newcommand{\EE}{{\EuScript{E}}}
	\newcommand{\HH}{{\EuScript{H}}}
	\newcommand{\JJ}{{\EuScript{J}}}
	\newcommand{\KK}{{\EuScript{K}}}
	\newcommand{\LL}{{\EuScript{L}}}
	\newcommand{\RR}{{\EuScript{R}}}
	\newcommand{\SS}{{\EuScript{S}}}

        \newcounter{sect}
        \newcommand{\sect}[1]{\vspace{3ex}\addtocounter{sect}{1}
                \begin{flushleft}
                {{\large\bf \arabic{sect}. {#1}}}
                \end{flushleft}
                \setcounter{thm}{0}
                \setcounter{equation}{0}
                \def\theequation{\arabic{sect}.\arabic{equation}}}
        \newcommand{\subsect}[2]{\begin{flushleft}
                {{\bf \arabic{sect}.{#1} {#2}}}
                \end{flushleft}}
        \newtheorem{thm}{Theorem}[sect]
        \newtheorem{prop}[thm]{Proposition}
        \newtheorem{lemma}[thm]{Lemma}
        \newtheorem{cor}[thm]{Corollary}

        \newcommand{\proof}[1]{\noindent {\em Proof.}$\quad$ {#1}
                        $\hfill\Box$
                        \vspace{1ex}}
        \newcommand{\be}{\begin{equation}}
        \newcommand{\ee}{\end{equation}}
        \newcommand{\bea}{\begin{eqnarray}}
        \newcommand{\eea}{\end{eqnarray}}
        \newcommand{\nno}{\nonumber \\}
        \newcommand{\sep}[1]{\!\!\!\! &{#1}& \!\!\!\! }
        \newcommand{\eq}{\sep{=}}
        \newcommand{\vc}{\sep{ }}

        \newcommand{\medwedge}{\mbox{\fontsize{12pt}{0pt}\selectfont $\wedge$}}
	\newcommand{\bra}{\langle}
	\newcommand{\ket}{\rangle}
	\newcommand{\ii}{\sqrt{-1}\,}
	\newcommand{\im}{{\mathrm Im}}
	\newcommand{\hf}{{\textstyle{\frac{1}{2}}}}
	\newcommand{\inv}[1]{{\frac{1}{{#1}}}}
	\newcommand{\sech}{\!\; {\mathrm sech}\,}  
	\newcommand{\transpose}[1]{{}^t\!#1}
	\newcommand{\tr}{{\mathrm Tr}}
	\newcommand{\diag}{{\mathrm diag}}
	\newcommand{\End}{{\mathrm End}}
	\newcommand{\siegel}{{\mathfrak{H}_n}}
	
	\newcommand{\abcd}{{\textstyle {A\;\;B \choose C\;\;D}}}
	\newcommand{\OmpOm}{\XI_{\Om'\Om}}
	\newcommand{\OmOmp}{\XI_{\Om\Om'}}
	\newcommand{\w}{{\Om}}

\begin{document}
  \begin{flushright}
math.SG/0409555
  \end{flushright}

  \vspace{7.5ex}

        \begin{center}
	{\LARGE\rm Geometric Quantization, Parallel Transport

\vspace{1ex}

	  and the Fourier Transform} \\

        \vspace{7ex}

        {\large William D.\ Kirwin\footnote{Current address:
Department of Mathematics and Statistics, Utah State University,
Logan, UT 84322-3900, USA; E-mail addresses: {\tt kirwin@math.usu.edu}} 
and Siye Wu\footnote{E-mail addresses: {\tt swu@euclid.colorado.edu}}}

        \vspace{1.75ex}
{\em Department of Mathematics, University of Colorado, Boulder, 
CO 80309-0395, USA}

        \end{center}
        \vspace{2ex}

\begin{abstract}
In quantum mechanics, the momentum space and position space wave
functions are related by the Fourier transform.
We investigate how the Fourier transform arises in the context of geometric
quantization.
We consider a Hilbert space bundle $\HH$ over the space $\JJ$ of compatible
complex structures on a symplectic vector space.
This bundle is equipped with a projectively flat connection.
We show that parallel transport along a geodesic in the bundle $\HH\to\JJ$
is a rescaled orthogonal projection or Bogoliubov transformation.
We then construct the kernel for the integral parallel transport operator.
Finally, by extending geodesics to the boundary (for which the metaplectic
correction is essential), we obtain the Segal-Bargmann and Fourier transforms
as parallel transport in suitable limits.\\

\end{abstract}


\sect{Introduction}
In quantum mechanics, the position and momentum space wave functions are 
related by the Fourier transform.
In this paper we examine how this relationship arises in the context of 
geometric quantization.
We provide a geometric interpretation of the Fourier transform as parallel
transport in a vector bundle of infinite rank.
In fact, this consideration leads to a smoothly parametrized family of
transforms which includes the Fourier transform, the Segal-Bargmann transform,
and the Bogoliubov transform. 

Quantization of a symplectic manifold $(M,\om)$ requires an Hermitian line
bundle $\ell\to M$ with a compatible connection such that the curvature is
$\frac\om\ii$.
$\ell$ is called a pre-quantum line bundle and it exists if and only if the 
de Rham class $[\frac{\om}{2\pi}]$ is integral.
The pre-quantum Hilbert space $\HH^0$ consists of sections of $\ell$ which are
square-integrable with respect to the Liouville volume form on $M$.  
As is well-known, $\HH^0$ is too large for the purpose of quantization.
The additional structure we need is an almost complex structure compatible
with $\om$.
The space $\JJ$ of such $J$ is connected and contractible.
Each $J\in\JJ$ defines a quantum Hilbert space $\HH_J$ of (square-integrable)
$J$-holomorphic sections of $\ell$.
They form a vector bundle $\HH$ of Hilbert spaces over $\JJ$, provided there
is no jump of $\dim\HH_J$ as $J$ varies.  

To compare $\HH_J$ with different $J$, we need a connection on $\HH\to\JJ$.
Given $J,J'\in\JJ$ and a path connecting them, parallel transport in $\HH$
is a unitary operator from $\HH_J$ to $\HH_{J'}$.
If the connection is projectively flat, the holonomy is $U(1)$, and parallel
transports along different paths from $J$ to $J'$ differ by at most a phase.
Since a quantum state is actually represented by a ray in the Hilbert space,
the ``physics'' obtained is thus independent of the choice of $J$.
Unfortunately such a connection does not exist in general [\ref{GM}].
We will henceforth restrict our attention to a symplectic vector space 
$(V,\om)$.
We also restrict $\JJ$ to the space of linear complex structures on $V$
compatible with $\om$.
In this case, a projectively flat connection on $\HH\to\JJ$ is constructed in
[\ref{ADPW}], as a finite dimensional model for studying Chern-Simons gauge
theory.

In this paper, we study parallel transport in the bundle $\HH$ along 
the geodesics in $\JJ$, when the symplectic manifold is a vector space.
The space $\JJ$ can be identified with the Siegel upper-half space, which
has a natural K\"ahler metric.  
We show that parallel transport along a geodesic in the Hilbert space bundle
is a rescaled orthogonal projection.
Hence parallel transport agrees with the Bogoliubov transformation in
[\ref{W1},\ref{W2}] and the intertwining operators in [\ref{Sa}] and [\ref{Ma}].
Part of the boundary of $\JJ$ (as a bounded domain) consists of real 
Langrangian subspaces $L$ of $(V,\om)$.
Each $L$ is a real polarization and also defines a quantum Hilbert space. 
By extending geodesics to the boundary (for which the metaplectic correction
is essential), we obtain the Segal-Bargmann and Fourier transforms as 
parallel transport in suitable limits.

The rest of the paper is organized as follows.
In Section~2, we recall the identification of $\JJ$ with the Siegel upper-half
space and describe the connection and the resulting geometry of the bundle of
quantum Hilbert spaces over $\JJ$.
We also incorporate the metaplectic correction.
In Section~3, we study parallel transport in $\HH$ along geodesics in $\JJ$.
The condition of parallel transport is a partial differential equation.
By the $Sp(V,\om)$ symmetry, it suffices to consider a special class of
geodesics so that the equation can be solved explicitly.
We then show that the parallel transport is actually a rescaled orthogonal
projection or the Bogoliubov transformation.
Hence the integral kernel for the equation of parallel transport is the
Bergman reproducing kernel, up to a positive factor.
We then show that by extending a geodesic in one direction to infinity,   
the parallel transport becomes the Segal-Bargmann transform.
Extending both ends of a geodesic to infinity, the parallel transport
converges to the Fourier transform.
Since a real Lagrangian space is on the boundary of $\JJ$, the quantum
Hilbert space associated to it is not inside the bundle $\HH\to\JJ$.  
We show that with the metaplectic correction, the limit holds in the
sense of almost everywhere convergence as sections over $V$.
Other ways to formulate the limit are also established.

Finally, we would like to mention some recent related work.
Let $K$ be a Lie group of compact type, that is, $K$ is locally isomorphic
to a compact Lie group.
The cotangent bundle $T^*K$ is naturally symplectic and, being diffeomorphic
to the complexification $K^\CO$, has a compatible complex structure.
In [\ref{Ha}], Hall constructed a generalized Segal-Bargmann transform between
the vertically polarized and K\"ahler polarized quantum Hilbert spaces.
The pairing is a unitary operator and a rescaled projection, as in 
[\ref{W1},\ref{W2}] for the flat case.
In [\ref{FMMN}], the authors study parallel transport in the quantum Hilbert
space bundle over a $1$-parameter family of K\"ahler polarizations on $T^*K$.
As in the flat case [\ref{ADPW}], the parallel transport equation is given by 
a holomorphic version of the heat operator, which also appeared in [\ref{Ha}].
It would be interesting to explore the projective flatness of the quantum
Hilbert space bundle over a larger class of K\"ahler polarizations on $T^*K$.

\sect{Geometry of the Hilbert Space Bundle}
\subsect{1}{Complex polarizations and the metaplectic correction}
Let $V$ be a real vector space of dimension $2n$ equipped with a constant
symplectic form $\om$ (i.e., a nondegenerate, closed 2-form).
There exist linear coordinates $\{x^i,y^i\}_{i=1,\dots,n}$ or
$\transpose{x}=(x^1\ \cdots\ x^n)$, $\transpose{y}=(y^1\ \cdots\ y^n)$
on $V$ such that
\[\om=\sum_{i=1}^n dx^i\wedge dy^i=\transpose{dx}\wedge dy.\]

A complex structure $J\in\End(V)$ is compatible with the symplectic form
$\om$ if $\om(J\,\cdot\,,J\,\cdot\,)=\om(\,\cdot\,,\,\cdot\,)$
and $\om(\,\cdot\,,J\,\cdot\,)>0$.
Given such a $J$, the complexification of $V$ decomposes as
\[V^\CO=V_J^{(1,0)}\oplus V_J^{(0,1)},\]
where $V_J^{(1,0)}=\{X\in V^\CO\,\vert\,JX=\ii X\}$
and $V_J^{(0,1)}=\overline{V_J^{(1,0)}}$.
Let $\JJ$ be the set of all compatible complex structures on $V$.
$\JJ$ can be identified, as follows, with the Siegel upper half-space
\[\siegel=\{\w\in M_n(\CO)\,\vert\,\transpose\Om=\Om,\ \im\Om>0\}
\subset\CO^{\inv{2}n(n+1)}.\]
We associate a compatible complex structure $J\in\JJ$ to a 
point $\w\in\siegel$ so that $V_J^{(1,0)}$ is the graph of $\w$.
Equivalently, the complex structure can be written in terms of 
$\Om=\Om_1+\ii\Om_2$ as
\[J={\w_1\w_2^{-1} \quad -\w_2-\w_1\w_2^{-1}\w_1 
\choose\w_2^{-1} \quad\quad\quad\quad -\w_2^{-1}\w_1}.\]
Thus $\JJ$ is identified with the positive Lagrangian Grassmannian. 
Real Lagrangian subspaces correspond to certain points on the boundary of 
$\siegel$.
For any $\w\in\siegel$, we choose the corresponding holomorphic coordinates
on $V$ as
\be\label{eqn:holocoord}
z_\w=(2\Om_2)^{-\inv{2}}(x-\bar\Om y).
\ee
The matrix factor $(2\Om_2)^{-\inv{2}}$ is chosen so that the symplectic form
is
\[\om=\ii\,\transpose{dz_\Om}\wedge d\bar z_\Om.\] 
We will drop the subscript $\Om$ when there is no danger of confusion.

There is a pre-quantum line bundle $\ell\to V$ with a connection whose
curvature is $\frac{\om}{\ii}$.
We use the symplectic potential $\tau=\hf\sum_{i=1}^nx^idy^i-y^idx^i$
to trivialize $\ell\to V$.
That is, the covariant derivative of a section $s\in\Gam(\ell)$ along $X$ is
\[\nabla_X s=X(s)-\ii(\iota_X\tau)s,\]
if $s$ is identified with a function on $V$. 
The pre-quantum Hilbert space $\HH^0$ consists of square-integrable
sections of $\ell$ with respect to the Liouville volume form 
$\veps_\om=\frac{\om^n}{(2\pi)^nn!}$.
Polarized sections of $\ell$ are those which are holomorphic, i.e.,
\[\nabla_{\bar z}\,s=0.\]
Using the complex coordinates \eqref{eqn:holocoord}, the covariant
derivatives in $\ell$ are
\be\label{eqn:covder}
   \nabla_z=\frac{\partial}{\partial z}-\inv{2}\bar z,\quad
   \nabla_{\bar z}=\frac{\partial}{\partial\bar z}+\inv{2}z.
\ee
Hence, a polarized section $\psi\in\HH_J$ can be written as
\[  \psi=\phi(z)\,e^{-\inv{2}|z|^2}  \]
for some entire function $\phi$.
Let $\HH_J\subset\HH^0$ denote the space of square integrable 
polarized sections with respect to the complex structure $J$.
This is the quantum Hilbert space.
We then have a quantum Hilbert space bundle $\HH\to\JJ$ with fiber
$\HH_J$ over $J$.
There is an Hermitian structure on this bundle given by
\be\label{eqn:innerprod}
\bra\psi_1,\psi_2\ket = \int_V\bar\psi_1\psi_2\,\veps_\om
\ee
for $\psi_1,\psi_2\in\HH_J$.
Here and below, when $J$ is parameterized by $\Om\in\siegel$, the subscript
$J$ can be replaced by $\Om$.
For example, we write $\HH_\Om=\HH_J$.

Since $\JJ$ is the positive Lagrangian Grassmannian, there is a natural
canonical bundle $\LL\to\JJ$ with fiber $V_J^{(1,0)}$ over $J$.
Let $\KK\to\JJ$ be the dual determinant bundle with fiber 
$\KK_J=\bigwedge^n(V_J^{(1,0)})^*$.
Since $\JJ$ is contractible, there is a unique (up to equivalence) square
root bundle $\sqrt\KK\to\JJ$ such that $\sqrt\KK\otimes\sqrt\KK=\KK$.
This square root bundle is known as the bundle of half-forms.
We define the corrected quantum Hilbert space bundle $\hat{\HH}\to\JJ$ as
$\hat\HH=\HH\otimes \sqrt\KK$.
The fiber $\hat\HH_J=\HH_J\otimes\sqrt\KK_J$ is called the corrected quantum
Hilbert space.
Including the bundle of half-forms is known as the metaplectic correction.

\subsect{2}{Symplectic and metaplectic group actions}
Given a vector space $V$ with a symplectic form $\om$, the symplectic group
$Sp(V,\om)$ is the set of linear transformations on $V$ preserving $\om$.
Upon choosing a set of linear symplectic coordinates 
$\{x^i,y^i\}_{i=1,\dots,n}$, the group $Sp(V,\om)$ is isomorphic to
\[Sp(2n,\RE)=\left\{\abcd\Big\vert \
\transpose{A}C=\transpose{C}A, \ \transpose{B}D=\transpose{D}B, \
\transpose{A}D-\transpose{C}B=I_n\right\}.\]
The group $Sp(V,\om)$ acts on the set $\JJ$ of compatible complex
structures by $g\colon J\mapsto gJg^{-1}$.
The corresponding action on positive complex Lagrangian subspaces is
$g\colon V_J^{(1,0)}\mapsto gV_J^{(1,0)}=V_{gJg^{-1}}^{(1,0)}$.
Identifying $\JJ$ with the Siegel upper half-space $\siegel$, 
the action of $Sp(V,\om)$ on $\JJ$ becomes the fractional linear
transformation on $\siegel$, i.e.,
\be
g=\abcd\colon\w\mapsto g\cdot\Om=(A\w+B)(C\w+D)^{-1}.
\ee
The following results, which will be used in the sequel, can be verified by
straightforward calculations.

\begin{lemma}\label{lemma:omomp}
Suppose $g=\abcd\in Sp(2n,\RE)$ and $\Om=g\cdot\Om_0,
\Om'=g\cdot\Om'_0\in\siegel$.
Put $\OmOmp=\frac{\Om-\bar\Om'}{2\ii}$.
Then\\
\vspace{-2ex}

1. 
\vspace{-10pt}
\be\label{eqn:omomp}
\Om-\bar\Om'=\transpose\overline{(C\Om'_0+D)}{}^{-1}(\Om_0-\bar\Om'_0) 
      (C\Om_0+D)^{-1}.
\ee
In particular,
\[ \im\Om=\transpose\overline{(C\Om_0+D)}{}^{-1}\im\Om_0\,(C\Om_0+D)^{-1}.\]
2. 
\vspace{-22.5pt}
\bea
\!\!\!\!\!(\Om_2^{-1}-\OmOmp^{-1})\Om_2=(\Om-\bar\Om')^{-1}(\bar\Om-\bar\Om')
  \eq(C\Om_0+D)(\Om_0-\bar\Om'_0)^{-1}(\bar\Om_0-\bar\Om'_0)
    \overline{(C\Om_0+D)}{}^{-1};                              \\
\!\!\!\!\!  \Om'_2({\Om'_2}^{-1}-\OmOmp^{-1})=(\Om-\Om')(\Om-\bar\Om')^{-1}
   \eq\transpose(C\Om'_0+D)^{-1}(\Om_0-\Om'_0)(\Om_0-\bar\Om'_0)^{-1}
    \transpose\overline{(C\Om'_0+D)}.
\eea
\end{lemma}

The action of $Sp(V,\om)$ on $V$ lifts to the pre-quantum line bundle $\ell$
preserving the connection.
Consequently, the group $Sp(V,\om)$ acts on the pre-quantum Hilbert space
$\HH^0$.
In fact, since the symplectic potential $\tau$ is preserved by $Sp(2n,\RE)$,
under the corresponding trivialization $\ell\cong V\times\CO$, the action 
of $g\in Sp(2n,\RE)$ is
\[ g\cdot(v,\zeta)=(gv,\zeta), \quad v\in V\cong\RE^{2n},\ \zeta\in\CO,  \]
and that on $s\in\HH^0\cong L^2(V)\otimes\CO$ is
\[  (g\cdot s)(v)=s(g^{-1}v), \quad v\in V.  \]
The action of $Sp(V,\om)$ lifts to the Hilbert space bundle $\HH\to\JJ$
covering the action on $\JJ$.
Since $Sp(V,\om)$ preserves the connection on $\ell$, the action
$g\colon\HH_J\to\HH_{gJg^{-1}}$ is a unitary isomorphism for any 
$g\in Sp(V,\om)$.

The symplectic group $Sp(V,\om)$ also acts on the vector bundle $\LL\to\JJ$
and hence on the line bundle $\KK\to\JJ$.
In fact the choice of coordinates \eqref{eqn:holocoord} provides a global
unitary section $\Om\mapsto d^nz_\Om$ of $\KK$.
The transformation of the complex coordinates
\be\label{eqn:ztrans}
g=\abcd\colon z_\Om\mapsto(g^{-1})^*z_\Om=\Om_2^{-\inv{2}}
\transpose\overline{(C\Om+D)}(g\cdot\Om)_2^{\inv{2}}z_{g\cdot\Om}=
\Om_2^{\inv{2}}(C\Om+D)^{-1}(g\cdot\Om)_2^{-\inv{2}}z_{g\cdot\Om},
\ee
where $(g\cdot\Om)_2=\im(g\cdot\Om)$, is unitary, and so is that of the
section $d^nz$
\be\label{eqn:gsqrt}
g=\abcd\colon d^nz_\Om\mapsto(g^{-1})^*d^nz_\Om
 =\frac{\det\overline{(C\Om+D)}}{|\det(C\Om+D)|}d^nz_{g\cdot\Om}.
\ee
This action does not lift to $\sqrt\KK$, but the double covering group of
$Sp(V,\om)$ does act on $\sqrt\KK$.
Since $Sp(V,\om)$ is connected and $\pi_1(Sp(V,\om))\cong\ZZ$,
there is a unique (up to isomorphism) connected double covering group
$Mp(V,\om)$, called the metaplectic group.
The double covering group of $Sp(2n,\RE)$ is denoted by $Mp(2n,\RE)$.
We have the following well-known result (see for example [\ref{RR}]):
\begin{prop} $Mp(V,\om)$ is isomorphic to the group whose elements are pairs
$(\sig,g)$, where $\sig$ is a bundle isomorphism of $\sqrt\KK\to\JJ$
covering the action of $g\in Sp(V,\om)$ on $\JJ$.
That is, we have a commutative diagram
\[\begin{array}{ccc}
\sqrt\KK & \stackrel{\sig}{\longrightarrow} & \sqrt\KK \\
\downarrow & & \downarrow \\
\JJ & \stackrel{g}{\longrightarrow} & \JJ
\end{array}\]
\end{prop}
Consequently, the metaplectic group $Mp(V,\om)$ acts on $\sqrt\KK$ and hence
acts on the corrected quantum Hilbert space bundle
$\hat\HH=\HH\otimes \sqrt\KK$.
Given $g=\abcd\in Sp(2n,\RE)$, the action of a lifted element in 
$Mp(2n,\RE)$ is
\be
\psi\otimes\sqrt{d^nz_\Om}\in\hat\HH_\Om\mapsto
\frac{\det\overline{(C\Om+D)}{}^{\inv{2}}}{|\det(C\Om+D)|^{\inv{2}}}
\,\psi\circ g^{-1}\otimes\sqrt{d^nz_{g\cdot\Om}}\in\hat\HH_{g\cdot\Om},
\ee
where the square root $\det\overline{(C\Om+D)}{}^{\inv{2}}$ depends on the
lift of $g$ to $Mp(2n,\RE)$.

\subsect{3}{Projectively flat and flat connections}
First, we describe a projectively flat connection on $\HH\to\JJ$ [\ref{ADPW}].
Combining this and the connection on $\sqrt\KK\to\JJ$, we obtain a flat 
connection on $\hat\HH\to\JJ$ [\ref{W2}, \S 10.2].

Since $\HH\to\JJ$ is a subbundle of the product bundle $\JJ\times\HH^0\to\JJ$,
the trivial connection on the latter projects to a connection on $\HH$.
This connection is [\ref{ADPW}]
\be\label{eqn:connH}
\nabla^\HH = \del+\inv{4}(\del J\,\om^{-1})^{ij}\,\nabla_{z^i}\,\nabla_{z^j},
\ee
where $\del$ is the exterior differential on $\JJ$.
The second term is a $1$-form on $\JJ$ valued in the set of skew-adjoint
operators on $\HH_J$.
Let $P_J=\hf(1-\ii J)\colon V^\CO\to V_J^{(1,0)}$ be the projection with
respect to the Hermitian form on $V^\CO$ defined by $\om$ and $J$.
Then the curvature of the above connection is [\ref{ADPW}]
\be\label{eqn:curvH}
F^\HH=-\inv{8}\tr(P_J\del J\wedge\del JP_J)\,{\mathrm id}_{\HH_J}.
\ee
So the connection is projectively flat [\ref{ADPW}].
Henceforth we omit the identity operator.

The connection described above blows up at the boundary of $\JJ$.
We will be interested in extending geodesics in $\JJ$ to the boundary.
In order to parallel transport along the extended geodesics in the next
section, we must employ the metaplectic correction.

The product bundle $V^\CO\times\JJ\to\JJ$ has an Hermitian structure defined
by $\om$ and $J$.
So a connection on the sub-bundle $\LL\to\JJ$ is given by the orthogonal
projection of the trivial connection.
Its curvature is
\be\label{eqn:curvL}
F^{\LL}=P_J\,\del(P_J\del P_J)=P_J\,\del P_J\wedge\del P_J\, P_J
=-\inv{4}P_J\,\del J\wedge\del J\,P_J.
\ee

\begin{prop} {\rm ([\ref{W2}, \S 10.2])} The induced connection on
the corrected quantum Hilbert space bundle $\hat\HH\to\JJ$ is flat.
\end{prop}
\proof{The connection on $\sqrt\KK$ is $F^{\sqrt\KK}=-\hf\tr\,F^{\LL}$.
So by \eqref{eqn:curvH} and \eqref{eqn:curvL}, 
\[F^\HH+F^{\sqrt\KK}=0.\]
The result was proved in [\ref{W2}, \S 10.2] using cocycles.}

The identification $\JJ=\siegel$ provides $\JJ$ with a convenient set of 
coordinates.
Using the variation of $J$,
\[\del J=\frac{\ii}{2}{\bar\Om\choose I_n}\Om_2^{-1}
  \del\Om\,\Om_2^{-1}(I_n,-\bar\Om)-\frac{\ii}{2}{\Om\choose I_n}
  \Om_2^{-1}\del\bar\Om\,\Om_2^{-1}(I_n,-\Om),\]
the connection \eqref{eqn:connH} becomes
\be\label{eqn:curvcalc}
\nabla^\HH = \del-\frac{\ii}{2}\transpose{\nabla_z}\Om_2^{-\inv{2}}
\del\bar\Om\,\Om_2^{-\inv{2}}\nabla_z.
\ee
The curvatures \eqref{eqn:curvH} is
\be\label{eqn:curvHK}
F^\HH=\inv{8}\tr(\Om_2^{-1}\del\Om\wedge\Om_2^{-1}\del\bar\Om).
\ee
The latter is proportional to the standard K\"ahler form on $\siegel$.
On the other hand, using the (unitary) global section $\sqrt{d^nz}$, the
connection on $\sqrt{\KK}$ is given by the $1$-form (for any $n\ge1$)
\be\label{eqn:connsqrt}
A^{\sqrt{\KK}}={\textstyle\frac{\ii}{4}}\tr(\Om_2^{-1}\del\Om_1).
\ee
Its curvature is the negative of \eqref{eqn:curvHK}.

The Hilbert space $\HH_J$ is the Fock space of a harmonic oscillator
with Hamiltonian $H=|z|^2$.
In the case $n=1$, the parameter $\Om$ is $\tau=\tau_1+\ii\tau_2$
in the upper half-plane.
A unitary basis for $\HH_J$ is 
$\{|k\ket=\frac{z^k}{\sqrt{k!}}e^{-\inv{2}|z|^2}\}_{k\in\NN}$.
The vector $|0\ket$ is the vacuum state and $|k\ket$ ($k\ge1$) are the
excited states.
Such a basis provides a global unitary frame for the bundle $\HH$.
Each $|k\ket$, regarded as a function of $\tau$ valued in $\HH^0$, has
the exterior derivative
\[ \del|k\ket=\frac{\ii}{4\tau_2}\left[
\left(k|k\ket-\sqrt{(k+1)(k+2)}|k+2\ket\right)\del\tau
+\left(k|k\ket-2\bar z\sqrt{k}|k-1\ket+\bar z^2|k\ket\right)\del\bar\tau
\right]. \]
The connection is given by an infinite skew-Hermitian matrix valued $1$-form
\be
A^\HH_{kl}=\bra k|\del|l\ket=\frac{\ii}{4\tau_2}
\left[\left(k\,\del_{kl}-\sqrt{k(k-1)}\,\del_{k,l+2}\right)\del\tau+
\left(l\,\del_{kl}-\sqrt{l(l-1)}\,\del_{k+2,l}\right)\del\bar\tau\right],
\ee
while the matrix of the curvature $2$-form is, as expected, 
\be
\bra k|F^\HH|l\ket=\frac{\del\tau\wedge\del\bar\tau}{8\tau_2^2}\,\del_{kl}.
\ee

\sect{Parallel Transport along the Geodesics}
\subsect{1}{Solutions to the equation of parallel transport}
The Siegel upper half-space $\siegel$ has a non-positively curved
K\"ahler metric
\[  ds^2=\tr(\Om_2^{-1}\,\del\Om\,\Om_2^{-1}\,\del\bar\Om),   \]
which is invariant under the action of $Sp(2n,\RE)$.
We study parallel transport in the bundles $\HH$ and $\hat{\HH}$ along the
geodesics in $\siegel$.
Let $\Om,\Om'\in\siegel$ represent $J,J'\in\JJ$, respectively.
The parallel transport in the bundle $\HH$ along the unique geodesic from
$\Om$ to $\Om'$ is a unitary operator, and so is that in $\hat\HH$.
We denote them by $U_{J'J}=U_{\Om'\Om}\colon\HH_J\to\HH_{J'}$ and
$\hat U_{J'J}=\hat U_{\Om'\Om}\colon\hat\HH_J\to\hat\HH_{J'}$, respectively.
The generating function for the basis of $\HH_\Om$ is a coherent state or
a principal vector [\ref{Ba}]
\be
c_\al(z_\Om)=\exp(\transpose\bar\al z_\Om-\textstyle{\inv{2}}|z_\Om|^2),
\ee
where $\al\in\CO^n$.
We wish to find $U_{\Om'\Om}\,c_\al\in\HH_{\Om'}$ and its metaplectic
correction.

For any diagonal matrix $\Lam=\diag[\lam_1,\dots,\lam_n]\ge0$, the curve
$\gam_\Lam\colon\RE\to\siegel$ defined by $\gam_\Lam(t)=\ii e^{2\Lam t}$
is a geodesic in $\siegel$.

\begin{lemma}\label{lemma:geodesics} {\rm ([\ref{Si}])}
For any geodesic $\gam\colon\RE\to\siegel$, there exist $g\in Sp(2n,\RE)$
and a diagonal matrix $\Lam\ge0$ such that $\gam=g\cdot\gam_\Lam$.
\end{lemma}

We first study parallel transport along the geodesic $\gam_\Lam$; the latter
determines a one-parameter family of complex structures $J_t$, whose complex
coordinates are $z_t=\inv{\sqrt{2}}(e^{-\Lam t}x+\ii e^{\Lam t}y)$.
The equation of parallel transport of a family of polarized sections
$\psi_t\in\Gam(\gam_\Lam^*\HH)$ is
\be\label{eqn:eqn}
(\partial_t-\hf\,\transpose{\nabla_{z_t}}\Lam\,\nabla_{z_t})\psi_t=0.
\ee

\begin{prop}\label{prop:sol}
The parallel transport of $c_\al(z_0)$ along $\gam_\Lam$ from $\Om_0=\ii I_n$
to $\Om_t=\ii e^{2\Lam t}$ is given by
\be\label{eqn:sol}
(U_{\Om_t\Om_0}\,c_\al)(z)=(\det\sech\Lam t)^\hf
\,\exp\left[\,\inv{2}\,{\bar\al\choose z}^{\hspace{-20pt}t\hspace{16.25pt}}
{\tanh\Lam t\quad\sech\Lam t\;\choose\sech\Lam t\;\;-\!\tanh\Lam t}
{\bar\al\choose z}-\inv{2}|z|^2\right].
\ee
\end{prop}

\proof{Since the connection $1$-form $-\hf\transpose\nabla_z\Lam\nabla_z$
is a sum of diagonal terms, we can assume $n=1$; the general case is similar.
We can also set $\lam_1=1$ by a rescaling of $t$.
Let $\psi_t$ be the parallel transport of $c_\al$.
Write $z=z_t$ and $\psi_t=\phi(t,z)e^{-\inv{2}|z|^2}$.
Then $\phi(t,z)$ is an entire function in $z$ (for each $t$) satisfying
\[\left(\frac{\partial}{\partial t}-\inv{2}\frac{\partial^2}{\partial z^2}
+\inv{2}z^2\right)\phi(t,z)=0,\quad\phi(0,z)=e^{\bar\al z}.\]
Here we have used \eqref{eqn:covder} and $\frac{d}{dt}z=-\bar z$.
Set $\phi(t,z)=e^{f(t,z)}$.
Then $f(t,z)$ satisfies
\[f_t-\hf f_{zz}-\hf f_z^2+\hf z^2=0,\quad f(0,z)=\bar\al z.\]
If we look for a solution of the form 
$f(t,z)=\hf p(t)z^2+q(t)z+\hf r(t)$, then $p$, $q$, $r$ satisfy a
set of ordinary differential equations
\[p_t=p^2-1,\quad q_t=pq,\quad r_t=p+q^2\]
with the initial conditions $p(0)=r(0)=0$, $q(0)=\bar\al$.
The solutions are $p(t)=-\tanh t$, $q(t)=\bar\al\sech t$, 
$r(t)=\ln\sech t+\bar\al^2\tanh t$, and hence
\[\phi(t,z)=\sqrt{\sech t}\,
\exp\left(\bar\al z\sech t+\hf(\bar\al^2-z^2)\tanh t\right).\]
The result follows from the uniqueness of parallel transport.}

Proposition~\ref{prop:sol} enables us to calculate the parallel transport 
of any basis vector in $\HH_{J_0}$.
In particular, the parallel transport of the vacuum is no longer the vacuum
in a new polarization; it is a linear combination of states with an even
number of excitations.
We list the parallel transport of a few states with small excitation numbers
in the case $n=1$, $\lam_1=1$:
\be
\left(\begin{array}{c}1\\ z_0\\ z_0^2\\ \vdots\end{array}\right)
\,e^{-\inv{2}|z_0|^2}\mapsto\sqrt{\sech t}
\left(\begin{array}{c}1\\ z_t\sech t \\ \!\!z_t^2\sech\!^2 t+\tanh t\!\!
\\ \vdots\end{array}\right)
\,e^{-\inv{2}z_t^2\tanh t-\inv{2}|z_t|^2}.
\ee

Next we study parallel transport in the half-form bundle $\sqrt\KK$.
As noted earlier, the complex coordinates corresponding to the point
$\gam_\Lam(t)\in\siegel$ are
$z_t=\inv{\sqrt{2}}(e^{-\Lam t}x+\ii e^{\Lam t}y)$.
As $t$ varies, the complex coordinates change by 
$\frac{d}{dt}z_t=-\Lam\bar z_t$, whose projection to 
$V_{\gam_\Lam(t)}^{(1,0)}$ is $0$.
Consequently, $dz_t$ is a parallel section of $\gam_\Lam^*\LL^*$
and $\sqrt{d^n z_t}$ is a parallel section of $\gam_\Lam^*\KK$.
The latter is also a consequence of \eqref{eqn:connsqrt}.
Hence the parallel transport of $c_\al\otimes\sqrt{d^n z_0}$ is 
$U_{\Om_t\Om_0}\,c_\al\otimes\sqrt{d^n z_t}$.

We now turn to parallel transport along a general geodesic.

\begin{thm}
Let $\Om,\Om'\in\siegel$ and let $\gam$ be the unique geodesic such that
$\gam(0)=\Om$ and $\gam(1)=\Om'$.
Then\\
1. The parallel transport of $c_\al(z_\Om)$ along $\gam$ is
{\small
\be
(U_{\Om'\Om}\,c_\al)(z'_{\Om'})=
 \frac{(\det\Om_2)^\inv{4}(\det\Om'_2)^\inv{4}}{|\det\OmOmp|^{\inv{2}}}
 \,\exp\!\!\left[\,\inv{2}\,
 {\bar\al\choose z'_{\Om'}}^{\hspace{-26pt}t\hspace{23pt}}
 {I_n-\Om_2^\inv{2}\OmOmp^{-1}\Om_2^\inv{2}\quad\;
 \Om_2^\inv{2}\OmOmp^{-1}{\Om'_2}^\inv{2}\quad
 \choose\quad{\Om'_2}^\inv{2}\OmOmp^{-1}\Om_2^\inv{2}
 \quad\;I_n-{\Om'_2}^\inv{2}\OmOmp^{-1}{\Om'_2}^\inv{2}} 
 {\bar\al\choose z'_{\Om'}} -\inv{2}|z'_{\Om'}|^2   \right].    \nno
\ee}
\vspace{-13.5pt}
\be\label{eqn:transcoh}
\ee
2. The parallel transport of $\sqrt{d^nz_\Om}$ along $\gam$ is
\be\label{eqn:transqrt}
\frac{(\det\OmpOm)^{\inv{2}}}{|\det\OmpOm|^{\inv{2}}}\sqrt{d^nz_{\Om'}}=
\frac{|\det\OmOmp|^{\inv{2}}}{(\det\OmOmp)^{\inv{2}}}\sqrt{d^nz_{\Om'}}.
\ee
3. The parallel transport of $c_\al(z_\Om)\otimes\sqrt{d^nz_\Om}$ along $\gam$
is the tensor product of \eqref{eqn:transcoh} and \eqref{eqn:transqrt}.
\end{thm}

\proof{Let $g={A\quad B\choose C\quad D}$ and $\Lam$ be given by 
Lemma~\ref{lemma:geodesics} such that $\gam=g\cdot\gam_\Lam$.
Then $\Om=g\cdot\Om_0,\Om'=g\cdot\Om'_0\in\siegel$, where $\Om_0=\ii I_n$
and $\Om'_0=\ii e^{2\Lam}I_n$.\\
1. We first map $c_\al(z_\Om)=\exp(\transpose\bar\al z_\Om-\hf|z_\Om|^2)$
in $\HH_\Om$ by $g^{-1}$ to 
$c_{\al_0}(z_0)=\exp(\transpose\bar\al_0 z_0-\hf|z_0|^2)$ in $\HH_{\Om_0}$.
By the unitarity of \eqref{eqn:ztrans}, we have $|z_0|^2=|z_\Om|^2$ and 
\[\al_0=\transpose\overline{(C\Om_0+D)}\,\Om_2^\inv{2}\al
     =(C\Om_0+D)^{-1}\Om_2^{-\inv{2}}\al.  \]
The parallel transport of $c_{\al_0}(z_0)$ in $\HH_{\Om_0}$ along 
$\gam_\Lam$ is $(U_{\Om'_0\Om_0}c_{\al_0})(z'_0)$ in $\HH_{\Om'_0}$ 
given by \eqref{eqn:sol}.
Since the connection is invariant under $Sp(2n,\RE)$, the action of $g$
on the latter is $(U_{\Om'\Om}\,c_\al)(z'_{\Om'})$ in $\HH_{\Om'}$.
Here 
\[z'_0=e^{-\Lam}\;\transpose\overline{(C\Om'_0+D)}\,{\Om'_2}^\inv{2}z'_{\Om'}
      =e^\Lam(C\Om'_0+D)^{-1}{\Om'_2}^{-\inv{2}}z'_{\Om'}. \]
Using these identities and Lemma~\ref{lemma:omomp}, we get 
\[ \det(\cosh\Lam t)=\det\XI_{\Om_0\Om'_0}\det(\im\Om'_0)^{-\inv{2}}
   =|\det\OmOmp|(\det\Om_2)^{-\inv{2}}(\det\Om'_2)^{-\inv{2}}, \]
\[ \transpose\bar\al_0\sech\Lam t\;z'_0
  =\transpose\bar\al\Om_2^\inv{2}(C\Om_0+D)\XI_{\Om_0\Om'_0}^{-1}
   \transpose\overline{(C\Om_0+D)}{\Om'_2}^\inv{2}z'_{\Om'}
  =\transpose\bar\al\Om_2^\inv{2}\OmOmp^{-1}{\Om'_2}^\inv{2}z'_{\Om'}, \]
\[ \transpose\bar\al_0\tanh\Lam t\;\bar\al_0
  =\transpose\bar\al\Om_2^\inv{2}(C\Om_0+D)(\Om_0-\bar\Om'_0)^{-1}
   (\bar\Om_0-\bar\Om'_0)\overline{(C\Om_0+D)}{}^{-1}\Om_2^{-\inv{2}}\bar\al
  =\transpose\bar\al(I_n-\Om_2^\inv{2}\OmOmp^{-1}\Om_2^\inv{2})\bar\al, \]
and
\[ \hspace{-7.5pt}-\transpose z'_0\tanh\Lam t\;z'_0
  =\transpose z'_{\Om'}{\Om'_2}^{-\inv{2}}\,\transpose(C\Om'_0+D)^{-1}
   (\Om_0-\Om'_0)(\Om_0-\bar\Om'_0)^{-1}\transpose\overline{(C\Om'_0+D)}
   {\Om'_2}^\inv{2}z'_{\Om'}
  =\transpose z'_{\Om'}(I_n-{\Om'_2}^\inv{2}\OmOmp^{-1}{\Om'_2}^\inv{2})
   z'_{\Om'}. \]
{}From these identities and from \eqref{eqn:sol}, the result follows.\\
2. Since the connection on $\sqrt{\KK}$ is invariant under $Mp(2n,\RE)$,
we again map $\sqrt{d^nz_\Om}$ by $g^{-1}$ to $\Om_0$, parallel transport
it along $\gam_\Lam$ to $\Om'_0$, and map the result to $\Om'$ by $g$.
By \eqref{eqn:gsqrt}, the phase accumulated in these steps is
\[ \left(\frac{\det(C\Om_0+D)\det\overline{(C\Om'_0+D)}}
   {|\det(C\Om_0+D)\det(C\Om'_0+D)|}\right)^{\!\!\inv{2}},  \]
which is equal to those in \eqref{eqn:transqrt} by taking the determinant of
\eqref{eqn:omomp}.}

\subsect{2}{Projections, Bogoliubov transformations and the
integral kernel of parallel transport}
Given any two compatible complex structures $J,J'\in\JJ$, there is an
orthogonal projection $P_{J'J}\colon\HH_J\to\HH_{J'}$ inside the pre-quantum
Hilbert space $\HH^0$.
In [\ref{W1}] and [\ref{W2}, \S 9.9], Woodhouse showed that up to a scalar
multiplication, this projection is a unitary operator, the Bogoliubov
transformation.
On the other hand, parallel transport in the bundle $\HH\to\JJ$ along
the geodesic from $J$ to $J'$ defines a manifestly unitary operator from
$\HH_J$ to $\HH_{J'}$.
The rescaled projection of the vacuum state calculated in [\ref{W2}, \S 9.9] 
coincides with \eqref{eqn:sol} when $\al=0$.
We show that this is true for all states.

\begin{thm}\label{thm:bog}
For any $J,J'\in\JJ$, the parallel transport $U_{J'J}$ in $\HH$
along the (unique) geodesic from $J$ to $J'$ is the map 
$\al(J,J')P_{J'J}\colon\HH_J\to\HH_{J'}$, where
\be\label{eqn:scalefactor}
\al(J,J')=\left(\det\mbox{{\large $\textstyle\frac{J+J'}{2}$}}\right)^\inv{4}
=\left(\det\ii(P_{J'}-\bar P_J)\right)^{\inv{4}}.
\ee
\end{thm}

\proof{The quantity $\al(J,J')$ defined in \eqref{eqn:scalefactor} in
invariant under $Sp(V,\om)$.
Since all the steps are equivariant under $Sp(V,\om)$, it suffices to
consider the parallel transport along a geodesic of the form $\gam_\Lam$.
Let $J_t$ be the complex structure corresponding to $\gam(t)$.
We want to show that there exists a function $\al(t)$ with $\al(0)=1$
such that for any $\psi'\in\HH_{J_0}$, $\psi_t\in\HH_{J_t}$, if $\{\psi_t\}$
is parallel along $\gam_\Lam$, then
\[\bra\psi',\psi_0\ket=\al(t)\bra\psi',\psi_t\ket\]
for any $t\in\RE$.
This is equivalent to the condition that the right hand side has vanishing
derivative with respect to $t$.
Again,  without loss of generality, we prove the case of $n=1$, $\lam_1=1$.
Since $z_t=\inv{\sqrt{2}}(e^{-t}x+\ii e^ty)$, we have
\[ \nabla_{z_t}=\sech t\,\nabla_{z_0}+\tanh t\,\nabla_{\bar z_t}.  \]
We also note that the formal adjoint of $\nabla_{z_0}$ is $-\nabla_{\bar z_0}$.
So
\bea
\textstyle{\frac{d}{dt}}\vc\hspace{-12.55pt}(\al(t)\bra\psi',\psi_t\ket)=
 \al'(t)\bra\psi',\psi_t\ket+\hf\al(t)\bra\psi',\nabla_{z_t}^2\psi\ket  \nno
 \eq\al'(t)\bra\psi',\psi_t\ket+\hf\al(t)\bra\psi',(\sech t\,\nabla_{z_0}+
 \tanh t\,\nabla_{\bar z_t})\nabla_{z_t}\psi_t\ket                \nno
\eq\al'(t)\bra\psi',\psi_t\ket
 +\hf\al(t)\sech t\bra-\nabla_{\bar z_0},\nabla_{z_t}\psi_t\ket
 +\hf\al(t)\tanh t\,\bra\psi',(-\ii\om(\partial_{\bar z_t},\partial_{z_t})
  +\nabla_{z_t}\nabla_{\bar z_t})\psi_t\ket                     \nno
\eq(\al'(t)-\hf\al(t)\tanh t)\bra\psi',\psi_t\ket,              \nonumber
\eea
which vanishes if we choose $\al(t)=\sqrt{\cosh t}$.
It is easy to verify $\al(J_0,J_t)=\al(t)$ using 
$J_t={\quad\quad -e^{2t} \choose e^{-2t}\quad\quad}$.
The second equality in \eqref{eqn:scalefactor} is because
$\frac{J+J'}{2}=\ii(P_{J'}-\bar P_J)$.}
 
\begin{cor}
Parallel transport in $\HH\to\JJ$ along the geodesic in $\siegel$
from $J$ to $J'$ coincides with the Bogoliubov transformation from $\HH_J$
to $\HH_{J'}$.
\end{cor}

\proof{The operator $\al(J,J')P_{J'J}$, including the scalar factor
\eqref{eqn:scalefactor}, coincides with the formula of the Bogoliubov
transformation in  [\ref{W1}] and [\ref{W2}, \S 9.9].}

There is a more direct explanation of the above result.
Given a complex structure $J$ and the corresponding holomorphic coordinates
\eqref{eqn:holocoord} on $V$, $\HH_J$ is the Fock space of the creation and
annihilation operators
\[ a_J^\dagger=z,\quad a_J=\nabla_z+\bar z. \]
As the complex structure changes along a geodesic, so do $a_J^\dagger$ and
$a_J$.
In fact, when $n=1$ and along $\gam(t)=\ii e^{2t}$, the parallel transports
of $a_0$ and $a_0^\dagger$ are
\[ \cosh t\,a_t+\sinh t\,a_t^\dagger\quad\mbox{and}\quad
   \sinh t\,a_t+\cosh t\,a_t^\dagger,\]
respectively, where $a_t=a_{J_t}$ and $a_t^\dagger=a_{J_t}^\dagger$. 
This deformation of the creation and annihilation operators, or the concept
of the vacuum and excitations, is the physics origin of the Bogoliubov
transformation.

For any $J\in\JJ$ represented by $\Om$, $\HH_J$ can be identified with the 
space of analytic function on $(V,J)$ with the measure 
$e^{-|z_\Om|^2}\eps_\om(z_\Om)$. 
The orthogonal projection onto $\HH_\Om$ is given by the Bergman kernel.
So we can express parallel transport by an integral kernel operator.

\begin{prop}\label{prop:bergman}
For any $\Om,\Om'\in\siegel$, the parallel transport $U_{\Om'\Om}$ in $\HH$
along the (unique) geodesic from $\Om$ to $\Om'$ is 
\be\label{eqn:bergman}
\phi(z_\Om)\,e^{-\inv{2}|z_\Om|^2}\mapsto
\frac{|\det\OmpOm|^\inv{2}}{(\det\Om_2)^\inv{4}(\det\Om'_2)^\inv{4}}
\;e^{-\inv{2}|z'_{\Om'}|^2}\int_V e^{\transpose{z'_{\Om'}}\bar z_{\Om'}
-\inv{2}|z_{\Om'}|^2-\inv{2}|z_\Om|^2}\phi(z_\Om)\,\eps_\om(z_\Om).
\ee
\end{prop}

\proof{The projection onto $\HH_{\Om'}$ is given by the Bergman kernel
\[ e^{\transpose{z'_{\Om'}}\bar z_{\Om'}
      -\inv{2}|z'_{\Om'}|^2-\inv{2}|z_{\Om'}|^2}.         \]
Using the facts $J{\Om \choose 1}=\ii{\Om \choose 1}$ and
$(1,-\Om)J=-\ii(1,-\Om)$, we get
\[ {1\,\;-\bar\Om\choose 1\,\;-\Om}\mbox{\Large $\textstyle\frac{J+J'}{2}$}
   {\Om'\quad\bar\Om'\choose 1\quad\;\;1}=\ii{\Om'-\bar\Om
   \quad\quad 0\;\quad\choose\quad\; 0\quad\quad\Om-\bar\Om'}. \]
Taking the determinant, we get
\[ \det\mbox{\large $\textstyle\frac{J+J'}{2}$}
   =\frac{|\det\OmpOm|^2}{\det\Om_2\det\Om'_2}, \]
from which the scalar factor in \eqref{eqn:bergman} follows.}

Up to a phase, \eqref{eqn:bergman} agrees with the unitary intertwining
operator from $\HH_\Om$ to $\HH_{\Om'}$ in [\ref{Sa}, \ref{Ma}].

A pairing can be defined on the half forms $\sqrt{d^nz_\Om}$ and 
$\sqrt{d^nz_{\Om'}}$ even though they come from different complex 
structures.\footnote{We recall that the pairing of $d^nz_\Om$ and 
$d^nz_{\Om'}$ is determined by 
$(-1)^{\frac{n(n-1)}{2}}\,\frac{d\bar z_{\Om'}\wedge dz_\Om}{(2\pi\ii)^n}
   =\bra d^nz_{\Om'},d^nz_\Om\ket\,\eps_\om $
and that of $\sqrt{d^nz_\Om}$ and $\sqrt{d^nz_{\Om'}}$ is
$ \bra\sqrt{d^nz_{\Om'}},\sqrt{d^nz_\Om}\,\ket
  =\sqrt{\bra d^nz_{\Om'},d^nz_\Om\ket}$.} 
A simple calculation yields
\be\label{eqn:hfpair}
\bra\sqrt{d^nz_{\Om'}},\sqrt{d^nz_\Om}\;\ket=
\frac{(\det\OmpOm)^\inv{2}}{(\det\Om_2)^\inv{4}(\det\Om'_2)^\inv{4}}.
\ee
Since both the scalar factor in \eqref{eqn:bergman} and the phase in
\eqref{eqn:transqrt} are absorbed in \eqref{eqn:hfpair}, we have recovered

\begin{cor}{\rm ([\ref{W2}, \S 10.2])}
Parallel transport from $\Om$ to $\Om'$ under the flat connection 
in $\hat\HH=\HH\otimes\sqrt{\KK}$ is given by
\be
 \hat U_{\Om'\Om}\colon\psi\otimes\sqrt{d^nz_\Om}\in\hat\HH_\Om\mapsto
 \bra\sqrt{d^nz_{\Om'}},\sqrt{d^nz_\Om}\;\ket P_{\Om'\Om}\psi
 \otimes\sqrt{d^nz_{\Om'}}\in\hat\HH_{\Om'}.
\ee
Alternatively, this map can be described by a pairing between $\hat\HH_\Om$
and $\hat\HH_{\Om'}$
\be
 \bra\psi'\otimes\sqrt{d^nz_{\Om'}},\psi\otimes\sqrt{d^nz_\Om}\,\ket
 =\bra\psi',\psi\ket\bra\sqrt{d^nz_{\Om'}},\sqrt{d^nz_\Om}\,\ket,
\ee
where $\bra\psi',\psi\ket$ is the inner product of $\psi\in\HH_\Om$ and
$\psi'\in\HH_{\Om'}$ in $\HH^0$.
\end{cor}

When $\Om'=\Om$, the above pairing is the inner product 
\eqref{eqn:innerprod}in $\hat\HH_\Om$.

We remark that if $\phi(z_\Om)\,e^{-\inv{2}|z_\Om|^2}$ in \eqref{eqn:bergman}
is $c_\al(z_\Om)$, the integration yields the same result as 
\eqref{eqn:transcoh}.
This gives another integral kernel of parallel transport.
The existence of two different kernels is because $\phi(z_\Om)$ is
restricted to be holomorphic. 

\begin{thm}
Under the assumptions of Proposition~\ref{prop:bergman}, the map
$U_{\Om'\Om}$ sends $\phi(z_\Om)\,e^{-\inv{2}|z_\Om|^2}$ to

{\small 
\be\hspace{-5pt}
\frac{(\det\Om_2)^\inv{4}(\det\Om'_2)^\inv{4}}{|\det\OmOmp|^{\inv{2}}}
 \;e^{-\inv{2}|z'_{\Om'}|^2}\int_V \phi(z_\Om)\exp\!\!\left[\,\inv{2}
 {{\bar z_\Om\choose z'_{\Om'}}^{\hspace{-26pt}t\hspace{23pt}}}
 {I_n-\Om_2^\inv{2}\OmOmp^{-1}\Om_2^\inv{2}\quad\;
 \Om_2^\inv{2}\OmOmp^{-1}{\Om'_2}^\inv{2}\quad
 \choose\quad{\Om'_2}^\inv{2}\OmOmp^{-1}\Om_2^\inv{2}
 \quad\;I_n-{\Om'_2}^\inv{2}\OmOmp^{-1}{\Om'_2}^\inv{2}} 
 {\bar z_\Om\choose z'_{\Om'}} -|z_\Om|^2   \right]
 \,\eps_\om(z_\Om).            \nno
\ee}

\vspace{-20pt}

\be\label{eqn:holoker}
\ee
\end{thm}

\proof{Since $\phi(z_\Om)e^{-\inv{2}|z_\Om|^2}$ is in $\HH_\Om$, we have an
estimate $|\phi(w)|\le C\,e^{\inv{2}|w|^2}$, where $C$ is its norm [\ref{Ba}].
By the reproducing property of the Bergman kernel, 
\[ \phi(z_\Om)\,e^{-\inv{2}|z_\Om|^2}
  =\int_V\phi(w)\,e^{-|w|^2}c_w(z_\Om)\,\eps_\om(w),   \]
which we substitute in \eqref{eqn:bergman}.
The integrand satisfies
\[ e^{-\inv{2}|z'_{\Om'}|^2}
\big\vert\,\phi(w)\,e^{-|w|^2-\transpose\bar wz_\Om-\inv{2}|z_\Om|^2
 +\transpose{z'_{\Om'}}\bar z_{\Om'}-\inv{2}|z_{\Om'}|^2}\big\vert
 \le C\,e^{-\inv{2}|w-z_\Om|^2-\inv{2}|z_{\Om'}-z'_{\Om'}|^2}.
\]
Hence the double integral in $w$ and $z_\Om$ is absolutely convergent.
Exchanging the order of the integration and performing the integral 
in $z_\Om$, we get
\[ \int_V\phi(w)\,e^{-|w|^2}(U_{\Om'\Om}\,c_w)(z'_{\Om'})\;\eps_\om(w), \]
which is \eqref{eqn:holoker} after relabelling the integration variable $w$
as $z_\Om$.}

\subsect{3}{Segal-Bargmann and Fourier transforms as parallel transport}
The set $\SS$ of real Lagrangian subspaces in $(V,\om)$ can be identified
with the Shilov boundary of $\JJ$.
For any $L\in\SS$, there is an Hermitian form on the space of sections of 
$\ell$ that are covariantly constant along $L$ by choosing a translation
invariant measure on $V/L$. 
The subspace $\HH_L$ of such sections that are $L^2$-integrable on $V/L$ 
is independent of the choice of the measure.
The bundles $\KK$ and $\sqrt\KK$ extend to $\SS$; the fiber of $\KK$
over $L\in\SS$ is $\KK_L=(\medwedge^n(V/L)^*)^\CO$, where $(V/L)^*$
is identified as the subspace of $V^*$ that annihilates $L$. 
Let $\hat\HH_L=\HH_L\otimes\sqrt\KK_L$ for any $L\in\SS$;
this is the quantum Hilbert space (with the metaplectic correction)
associated to the real polarization $L$.
The action of $Sp(V,\om)$ on $\SS$ lifts to that of $Mp(V,\om)$ on the
bundles $\sqrt\KK$ and $\hat\HH$.
There is a canonical Hermitian form on $\hat\HH$.
Given $\psi_1,\psi_2\in\HH_L$ and $\sqrt\nu\in\sqrt\KK_L$ ($L\in\SS$),
we have  
\be\label{eqn:prodL}
\bra\psi_1\otimes\sqrt\nu,\psi_2\otimes\sqrt\nu\,\ket=
\int_{V/L}\bar\psi_1\psi_2\,\textstyle\frac{|\nu|}{(2\pi)^{\frac{n}{2}}},
\ee
where $|\nu|$ is a density on $V/L\cong\RE^n$ determined by $\sqrt\nu$.

The $Mp(V,\om)$-invariant pairing on $\sqrt\KK\to\JJ$ between different 
fibers also extends.
Pairings are defined between $\sqrt\KK_J$, $\sqrt\KK_L$ and between 
$\sqrt\KK_L$, $\sqrt\KK_{L'}$ for $J\in\JJ$ and $L,L'\in\SS$ such that $L$
and $L'$ are transverse.
For example, if $L_-=\{x=0\},L_+=\{y=0\}\in\SS$ in the symplectic 
coordinates $(x,y)$ and $\Om\in\siegel$, we have
\be
   \bra\sqrt{d^nz_\Om},\sqrt{d^nx}\,\ket=\det((2\Om_2)^{-\inv{2}}
   \textstyle{\frac\Om\ii}),
   \quad \bra\sqrt{d^ny},\sqrt{d^nx}\,\ket=\ii^{\frac{n}{2}}.
\ee
For any $J\in\JJ$ corresponding to $\Om\in\siegel$, let $\RR_J$ be the
subspace of $\psi\in\HH_J$ such that
\[ |\psi(z_\Om)|\le\frac{C}{(1+|z_\Om|^2)^{n+\al}} \]
for some $C\ge0$ and $\al>0$; such a $\psi$ is $L^1$ on $V$.
Let $\hat\RR_J=\RR_J\otimes\sqrt\KK_J$.
There is a pairing
\be\label{eqn:pairJL}
\bra\psi\otimes\sqrt\nu,\psi'\otimes\sqrt{\nu'}\,\ket=
   \bra\sqrt\nu,\sqrt{\nu'}\,\ket\,\int_V\bar\psi\psi'\;\eps_\om  
\ee
between any $\psi\otimes\sqrt\nu\in\hat\RR_J$ and 
$\psi'\otimes\sqrt{\nu'}\in\hat\HH_L$; the integral in \eqref{eqn:pairJL}
is absolutely convergent.
The corresponding operator $\hat B_{JL}\colon\hat\HH_L\to\hat\HH_J$ is unitary
and intertwines with the $Mp(V,\om)$-action [\ref{Ba}, \ref{Ma}, \ref{W2}].
If $L=L_-$ and if $J$ is parameterized by $\Om\in\siegel$, the operator and
its inverse are, respectively,
{\small
\bea
\hat B_{\Om L}\colon\vc \phi(x)\,e^{\frac{\ii}{2}\transpose xy}
\otimes\sqrt{d^nx}\;\longmapsto\;
\frac{(\det2\Om_2)^\inv{4}}{\overline{(\det\frac\Om\ii)}{}^\inv{2}}
 \;e^{-\inv{2}|z'_\Om|^2}\!\int_{V/L_-}\!\phi(x)                    \nno 
\vc\exp\!\!\left[\,\inv{2}{z'_\Om\choose x}^{\hspace{-22pt}t\hspace{21pt}} \!
 {I_n-(2\Om_2)^\inv{2}\overline{(\frac\Om\ii)}{}^{-1}(2\Om_2)^\inv{2}\quad
 (2\Om_2)^\inv{2}\overline{(\frac\Om\ii)}{}^{-1} \choose\quad
 \overline{(\frac\Om\ii)}{}^{-1}(2\Om_2)^\inv{2}\quad\quad\quad\quad
 \;-\overline{(\frac\Om\ii)}{}^{-1}\!}{z'_\Om\choose x} \right]
\frac{|d^nx|}{(2\pi)^\frac{n}{2}}\otimes\sqrt{d^nz'_\Om},          \nonumber
\eea}

\vspace{-20pt}

\be\label{eqn:Btrans}
\ee
and, if $\phi(z_\Om)\,e^{-\inv{2}|z_\Om|^2}$ is in $\RR_\Om$,
{\small
\bea
\hat B_{\Om L}^{-1}\colon\vc\phi(z_\Om)\,e^{-\inv{2}|z_\Om|^2}
\otimes\sqrt{d^nz_\Om}\;\longmapsto\;
\frac{(\det2\Om_2)^\inv{4}}{(\det\frac\Om\ii)^\inv{2}}
 \,e^{\frac{\ii}{2}\transpose x'y'}\!\int_V \phi(z_\Om)\;e^{-|z_\Om|^2}  \nno 
\vc\exp\!\!\left[\,\inv{2}{\bar z_\Om\choose x'}^{\hspace{-22pt}t\hspace{21pt}}
 \!{I_n-(2\Om_2)^\inv{2}(\frac\Om\ii)^{-1}(2\Om_2)^\inv{2}\quad
 (2\Om_2)^\inv{2}(\frac\Om\ii)^{-1} \choose\quad
 (\frac\Om\ii)^{-1}(2\Om_2)^\inv{2}\quad\quad\quad\quad
 \;-(\frac\Om\ii)^{-1}\!}{\bar z_\Om\choose x'}\right]\eps_\om(z_\Om)
\otimes\sqrt{d^nx'}.                            \nonumber
\eea}
\nobreak
\vspace{-37pt}
\be\label{eqn:invB}
\vspace{10pt}
\ee
When $\Om=\ii I_n$, they are the usual Segal-Bargmann transform and
its inverse [\ref{Ba}].

For any pair of Lagrangian subspaces $L,L'\in\SS$ that are transverse,
there exists a Fourier transform operator
$\hat F_{L'L}\colon\hat\HH_L\to\hat\HH_{L'}$
that intertwines with the action of $Mp(V,\om)$ [\ref{LV}].
In particular, we have
\be\label{eqn:Ftrans}
\hat F_{L_+L_-}\colon\phi(x)\,e^{\frac{\ii}{2}\transpose xy}\otimes\sqrt{d^nx}
\;\longmapsto\;\ii^\frac{n}{2}\left(\int_{V/L_-}\!\phi(x)
\,e^{\ii\transpose xy'}\textstyle{\frac{|d^nx|}{(2\pi)^\frac{n}{2}}}\right)
e^{-\frac{\ii}{2}\transpose x'y'}\otimes\sqrt{d^ny'},
\ee
where the integral in the bracket is the usual Fourier transform 
$\tilde\phi(y')$.
Strictly speaking, \eqref{eqn:Ftrans} is valid only on the dense subspace
of $\psi\otimes\sqrt{\nu}\in\hat\HH_L$ such that $|\psi|$ is $L^1$ on $V/L$;
the operator then extends continuously to $\hat\HH_L$.

\begin{prop}
1. Let $J\in\JJ$ and let $L,L'\in\SS$ be transverse to each other.
Then for any $\hat\psi\in\hat\RR_J$ and $\hat\psi'\in\hat\RR_L$,
\be\label{eqn:limitR}
\lim_{J'\to L}\hat U_{J'J}\hat\psi=\hat B_{JL}^{-1}\hat\psi,  \quad
\lim_{J'\to L'}\hat B_{J'L}\hat\psi'=\hat F_{L'L}\hat\psi';
\ee
here the limit is pointwise in $V$ as $J'\to L$ or $L'$ from inside $\JJ$.\\
2. For any $J,J'\in\JJ$ and $L,L',L''\in\SS$ that are mutually transverse,
we have
\be\label{eqn:identify}
\hat B_{J'L}=\hat U_{J'J}\circ\hat B_{JL},\quad
\hat F_{L'L}=\hat B_{JL'}^{-1}\circ\hat B_{JL},\quad
\hat F_{L''L}=\hat F_{L''L'}\circ\hat F_{L'L}.
\ee
\end{prop}

\proof{1. Let $J,J'$ be parameterized by $\Om,\Om'\in\siegel$.
Without loss of generality, assume $L=L_-$.
Then the limit $J'\to L$ is $\Om'\to0$ with $\Om'_2>0$.
If $\hat\psi(z_\Om)=\phi(z_\Om)\,e^{-\inv{2}|z_\Om|^2}\otimes\sqrt{d^nz_\Om}$,
then $(\hat U_{\Om'\Om}\hat\psi)(z'_{\Om'})$ is the tensor product of
\eqref{eqn:holoker} and \eqref{eqn:transqrt}.
As $\Om'\to0$, $(\det2\Om'_2)^\inv{4}\sqrt{d^nz'_{\Om'}}\to\sqrt{d^nx'}$
and the integrand in \eqref{eqn:holoker} goes to that in \eqref{eqn:invB}.
Since the latter is absolutely integrable, the limit commutes with the
integration and thus the first limit in \eqref{eqn:limitR} follows.
We remark here that the scalar factor $(\det2\Om'_2)^\inv{4}$ that goes to
zero in the limit is absorbed by the half-form $\sqrt{d^nz'_{\Om'}}$.
The proof of the second limit is similar.\\
2. Since the connection on the bundle $\hat\HH\to\JJ$ is flat and
since $\hat U_{J'J}$ is the parallel transport from $J$ to $J'$,
we have $\hat U_{J''J'}\circ\hat U_{J'J}=\hat U_{J''J}$.
Using $\hat U_{J'J}\colon\hat\RR_J\to\hat\RR_{J'}$ and taking $J''\to L$,
we get $\hat B_{J'L}^{-1}\circ\hat U_{J'J}=\hat B_{JL}^{-1}$ on $\hat\RR_J$,
and hence on $\hat\HH_J$.
The proof of the other two identities are similar.}

We thus proved that, as $J'\to L\in\SS$, parallel transport of 
$\hat\psi\in\hat\RR_J$ from $J$ to $J'$ goes to $\hat B_{JL}^{-1}\hat\psi$.
Since the latter is not $L^2$ on $V$ and its norm is defined instead
by \eqref{eqn:prodL}, it is not obvious why the ``operator''  
$\lim_{J'\to L}\hat U_{J'J}$ is continuous on $\hat\HH_J$ or why
its image is contained in $\hat\HH_L$.
We now take the limit $J'\to L$ as $J'$ follows the path of a geodesic.

\begin{lemma} 
1. Let $\Lam\ge0$ be a diagonal matrix and $\gam=g\cdot\gam_\Lam$, a
geodesic in $\JJ$.
Then $\lim_{t\to\pm\infty}\gam(t)$ are real Lagrangian subspaces 
if and only if $\Lam>0$.\\
2. For any $J\in\JJ$ and $L\in\SS$, there is a geodesic $\gam$ in $\JJ$
such that $\gam(0)=J$, $\lim_{t\to-\infty}\gam(t)=L$.\\ 
3. A pair of real Lagrangian subspaces $L,L'$ are transverse
if and only if there is a geodesic $\gam$ in $\JJ$ such that 
$\lim_{t\to-\infty}\gam(t)=L$, $\lim_{t\to+\infty}\gam(t)=L'$.
\end{lemma}

\proof{1. Using the identification of $\JJ$ and $\siegel$, 
$\gam_\Lam(-\infty)=0$ and $\gam_\Lam(+\infty)=+\ii\infty\,I_n$
if and only if $\Lam>0$, in which case they are real Lagrangian subspaces
$L_-$ and $L_+$, respectively.
The result follows from the transitivity of the $Sp(V,\om)$ action on $\SS$.\\
2. Without loss of generality, assume $J$ is represented by $\Om=\ii I_n$.
Then for any diagonal $\Lam>0$, $\gam_\Lam(0)=J$ and 
$\lim_{t\to-\infty}\gam_\Lam(t)=L_-$.
The isotropic subgroup of $J$ in $Sp(V,\om)$ is isomorphic to $U(n)$ and
acts transitively on $\SS$.
Hence the result.\\
3. Let $\gam=g\cdot\gam_\Lam$ ($\Lam>0$) be the geodesic such that the
limits hold.
Then $L=g\,L_-$ and $L'=g\,L_+$.
$L$, $L'$ are transverse since $L_-$, $L_+$ are. 
Conversely, if $L$, $L'$ are transverse, then there exists $g\in Sp(V,\om)$
such that $L=g\,L_-$, $L'=g\,L_+$.
The geodesic $\gam=g\cdot\gam_\Lam$ for any $\Lam>0$ satisfies the
requirement.}

\begin{prop}
Let $\gam$ be a geodesic in $\JJ$ such that $\gam(0)=J$ and 
$\gam(-\infty)=L,\gam(+\infty)=L'\in\SS$.
Then for any $\hat\psi\in\hat\HH_J$, we have
\be\label{eqn:lim+-}
\lim_{t\to-\infty}\hat U_{\gam(t)J}\hat\psi=\hat B_{JL}^{-1}\hat\psi,  \quad
\lim_{t\to+\infty}\hat U_{\gam(t)J}\hat\psi=\hat F_{L'L}
\lim_{t\to-\infty}\hat U_{\gam(t)J}\hat\psi
\ee
almost everywhere on $V$.
\end{prop}

\proof{Without loss of generality, we assume $\gam=\gam_\Lam$ ($\Lam>0$).
Then $J$ is given by $\Om=\ii I_n$ and $L=L_-$, $L'=L_+$, while at $\gam(t)$,
$z_t=\inv{\sqrt{2}}(e^{-\Lam t}x+\ii e^{\Lam t}y)$.
Let $\pi_\pm\colon V\to V/L_\pm$ be the projections.
Let $(\hat B_{JL}^{-1}\hat\psi)(x,y)=\phi(x)\,e^{\frac{\ii}{2}\transpose xy}
\otimes\sqrt{d^nx}$.
Using \eqref{eqn:limitR} and \eqref{eqn:Btrans}, we get
\bea
{}\hspace{-30pt}\vc{}\hspace{-15pt}(\hat U_{\gam(t)J}\hat\psi)(x',y')
 =(\hat B_{\gam(t)L_-}\hat B_{JL_-}^{-1}\hat\psi)(x',y')                 \nno
{}\hspace{-30pt}\eq\det e^{-\Lam t}\left(\int_{V/L_-}\!\!\phi(x)\,
   e^{-\inv{2}\transpose(x-x')e^{-2\Lam t}(x-x')+\ii\transpose(x-x')y'}
   \textstyle{\frac{|d^nx|}{(2\pi)^\frac{n}{2}}}\!\right)
   e^{\frac{\ii}{2}\transpose x'y'}
   \otimes(\det\sqrt 2e^{\Lam t})^\inv{2}\sqrt{d^nz'_t} \label{eqn:limit-} \\
{}\hspace{-30pt}
\eq\left(\int_{V/L_-}\!\!\phi(x)\,e^{-\inv{2}\transpose(x-x')e^{-2\Lam t}(x-x')
   +\ii\transpose xy'}\textstyle{\frac{|d^nx|}{(2\pi)^\frac{n}{2}}}\!\right)
   e^{-\frac{\ii}{2}\transpose x'y'}
   \otimes(\det\sqrt 2e^{-\Lam t})^\inv{2}\sqrt{d^nz'_t}. \label{eqn:limit+}
\eea
As $t\to-\infty$, \eqref{eqn:limit-} goes to 
$\phi(x')\,e^{\frac{\ii}{2}\transpose x'y'}\otimes\sqrt{d^nx}$ pointwise
on $\pi_-^{-1}(E_\phi)$, where $E_\phi$ is the Lebesgue set of $\phi$
(see for example [\ref{SW}, Theorem~I.1.25] or [\ref{Fo}, Theorem~8.62]).
Again, the scalar factor $(\det\sqrt 2e^{\Lam t})^\inv{2}$ that vanishes
in the limit is absorbed by $\sqrt{d^nz'_t}$.
As $t\to+\infty$, \eqref{eqn:limit+} goes to 
$\ii^\frac{n}{2}\tilde\phi(y')\otimes\sqrt{d^ny'}$ pointwise on 
$\pi_+^{-1}(E_{\tilde\phi})$ (see for example [\ref{Fo}, Theorem~8.31(c)]);
this also follows from the $t\to-\infty$ limit by making an $Sp(V,\om)$
transformation that fixes $J$ and exchanges $L_+$ and $L_-$.
It is well known that the Lebesgue set of an $L^2$ function is the complement
of a measure-zero subset (see for example [\ref{Fo}, Theorem~3.20] or
[\ref{SW}, pp.\ 12-13]).}

We remark that since the elements in $\hat\HH_L$ ($L\in\SS$) are defined 
up to a set of measure-zero, the limits in \eqref{eqn:lim+-} are the best 
possible results for pointwise convergence.
The integral in \eqref{eqn:limit-}, being the convolution of $\phi$ and the
heat kernel, goes to $\phi(x')$ when $t\to-\infty$ as tempered distributions
on $V/L_-$ (see for example [\ref{Ta}, Proposition~3.5.1] or 
[\ref{Fo}, Corollary~8.46]).
In the same sense, the integral in \eqref{eqn:limit+} goes to $\tilde\phi(y')$
when $t\to+\infty$.
Hence the limits in \eqref{eqn:lim+-} hold as tempered distributions on $V$,
with the given trivialization of $\ell$.

Finally, we consider the limit in $L^2$-spaces.
$\SS$ is part of the topological boundary of $\JJ$ as a bounded domain.
We define a topology on the disjoint union $\EE$ of all $\hat\HH_L$
($L\in\SS$) and the total space of $\hat\HH\to\JJ$.
There is a bijection from $\EE$ to $(\JJ\sqcup\SS)\times\hat\HH_{J_0}$
if we fix any $J_0\in\JJ$.
The maps from $\hat\HH_J$ ($J\in\JJ$) and $\hat\HH_L$ ($L\in\SS$) to
$\hat\HH_{J_0}$ are $\hat U_{J_0J}$ and $\hat B_{J_0L}^{-1}$, respectively.
The space $\EE$ thus inherits the product topology on
$(\JJ\sqcup\SS)\times\hat\HH_{J_0}$.

\begin{cor}
Let $J\in\JJ$ and let $L,L'\in\SS$ be a transverse pair.
Then for any $\hat\psi\in\hat\HH_J$, in the above topology on $\EE$, we
have
\[ \lim_{J'\to L}\hat U_{J'J}\hat\psi=\hat B_{JL}^{-1}\hat\psi,\quad
   \lim_{J'\to L'}\hat U_{J'J}\hat\psi=
   \hat F_{L'L}\lim_{J'\to L}\hat U_{J'J}\hat\psi.  \]
\end{cor}

\proof{The limits follow directly from \eqref{eqn:identify}.}

\medskip

\noindent
{\small {\bf Acknowledgements.}
We are grateful to Arlan Ramsay for many helpful conversations and
suggestions -- in particular regarding the limits of the parallel transport
operator.  We would also like to thank Wicharn Lewkeeratiyutkul for
bringing to our attention the paper [\ref{FMMN}].}

\bigskip

        \newcommand{\athr}[2]{{#1}.\ {#2}}
        \newcommand{\au}[2]{\athr{{#1}}{{#2}},}
        \newcommand{\an}[2]{\athr{{#1}}{{#2}} and}
        \newcommand{\jr}[6]{{#1}, {\it {#2}} {#3}\ ({#4}) {#5}-{#6}}
        \newcommand{\pr}[3]{{#1}, {#2} ({#3})}
        \newcommand{\bk}[4]{{\it {#1}}, {#2}, ({#3}, {#4})}
        \newcommand{\cf}[8]{{\it {#1}}, {#2}, {#5},
                 {#6}, ({#7}, {#8}), pp.\ {#3}-{#4}}
        \vspace{3ex}
        \begin{flushleft}
{\bf References}
        \end{flushleft}
{\small
        \baselineskip=13pt
        \begin{enumerate}

\item\label{ADPW}
\au{S}{Axelrod} \an{S}{Della Pietra} \au{E}{Witten}
\jr{Geometric quantization of Chern-Simons gauge theory}
{J.\ Diff.\ Geom.}{33}{1991}{787}{902}

\item\label{Ba}
\au{V}{Bargmann}
\jr{On a Hilbert space of analytic functions and an associated integral
transform}{Comm.\ Pure Appl.\ Math.}{14}{1961}{187}{214}

\item\label{Fo}
\au{G.\ B}{Folland}
\bk{Real analysis.\ Modern techniques and their applications}
{John Wiley \& Sons}{New York}{1984} 

\item\label{FMMN}
\au{C}{Florentino} \au{P}{Matias} \an{J}{Mourao} \au{J.\ P}{Nunes}
{\it Geometric quantization, complex structures and the coherent state 
transform}, {\tt math.DG/0402313}

\item\label{GM}
\an{V.\ L}{Ginzburg} \au{R}{Montgomery}
\cf{Geometric quantization and no-go theorems}
{Poisson geometry (Warsaw, 1998), Banach Center Publ., 51}{69}{77}
{eds.\ \an{J}{Grabowski} \athr{P}{Urba\'nski}}{Polish Acad.\ Sci.}
{Warsaw}{2000} 

\item\label{Ha}
\au{B.\ C}{Hall}
\jr{Geometric quantization and the generalized Segal-Bargmann transform for
Lie groups of compact type}{Comm.\ Math.\ Phys.}{226}{2002}{233}{268}

\item\label{LV}
\an{G}{Lion} \au{M}{Vergne}
\bk{The Weil representation, Maslov index and theta series,
Prog.\ in Math.\ 6}{Birkh\"auser}{Boston, MA}{1980}, Part I

\item\label{Ma}
\au{B}{Magneron}
\jr{Spineurs symplectiques purs et indice de Maslov de plan Lagrangiens 
positifs}{J.\ Funct.\ Anal.}{59}{1984}{90}{122}

\item\label{RR}
\an{P.\ L}{Robinson} \au{J.\ H}{Rawnsley}
\bk{The metaplectic representation, Mp$^c$ structures and geometric
quantization, {\rm Mem.\ Amer.\ Math.\ Soc.\ Vol.81, No.410}}
{Amer.\ Math.\ Soc.}{Providence, RI}{1989}

\item\label{Sa}
\au{I}{Satake}
\jr{On unitary representations of a certain group extension 
{\rm (in Japanese)}}{Sugaku}{21}{1969}{241}{253};
\cf{Fock representations and theta-functions}
{Advances in the theory of Riemann surfaces}{393}{405}
{\athr{L.\ V}{Ahlfors} et al eds.}{Princeton Univ.\ Press}{Princeton, NJ}{1971}

\item\label{Si}
\au{C.\ L}{Siegel}
\jr{Symplectic geometry}{Amer.\ J.\ Math.}{65}{1943}{1}{86}

\item\label{SW}
\an{E.\ M}{Stein} \au{G}{Weiss}
\bk{Introduction to Fourier analysis on Euclidean spaces}
{Princeton Univ.\ Press}{Princeton, NJ}{1971}

\item\label{Ta}
\au{M.\ E}{Taylor}
\bk{Partial differential equations I.\ basic theory}
{Springer-Verlag}{New York}{1996}

\item\label{W1}
\au{N.\ M.\ J}{Woodhouse}
\jr{Geometric quantization and the Bogoliubov transformation}
{Proc.\ Royal Soc.\ London A}{378}{1981}{119}{139}

\item\label{W2}
\au{N.\ M.\ J}{Woodhouse}
\bk{Geometric Quantization (2nd ed.)}
{Oxford Univ.\ Press}{New York}{1992}

\end{enumerate}}
\end{document}